\documentclass[11pt]{article}
\usepackage{amssymb,amsmath,amsthm}
\usepackage{epsfig}
\newtheorem{Theorem}{\sc Theorem}[section]
\newtheorem{Lemma}[Theorem]{\sc Lemma}
\newtheorem{Proposition}[Theorem]{\sc Proposition}
\newtheorem{Corollary}[Theorem]{\sc Corollary}
\newtheorem{Definition}[Theorem]{\sc Definition}
\newtheorem{Example}[Theorem]{\sc Example}
\newtheorem{Remark}[Theorem]{\sc Remark}

\def\hpic #1 #2 {\mbox{$\begin{array}[c]{l} 
\epsfig{file=#1,height=#2}\end{array}$}}
\def\wpic #1 #2 {\mbox{$\begin{array}[c]{l} 
\epsfig{file=#1,width=#2}\end{array}$}}

\def\C{\mathbb C}

\def\R{\mathbb R}

\def\CE{{\cal {E}}}
\def\CF{{\cal {F}}}

\def\CH{{\cal {H}}}

\def\CK{{\cal {K}}}
\def\CL{{\cal {L}}}

\def\CP{{\cal {P}}}

\def\CR{{\cal {R}}}

\def\be{\begin{equation}}
\def\ee{\end{equation}}
\def\bt{\begin{Theorem}}
\def\et{\end{Theorem}}
\def\bi{\begin{itemize}}
\def\ei{\end{itemize}}
\def\bea{\begin{eqnarray}}
\def\eea{\end{eqnarray}}
\def\ba{\begin{array}}
\def\ea{\end{array}}
\def\beast{\begin{eqnarray*}}
\def\eeast{\end{eqnarray*}}
\def\ben{\begin{enumerate}}
\def\een{\end{enumerate}}

\def\bi{\bibitem}

\def\rar{\rightarrow}
\def\Rar{\Rightarrow}

\def\i{{\bf i}}
\def\j{{\bf j}}

\def\u{{\overline{u}}}

\begin{document}

\begin{center}
{\Large {\bf On extendability of endomorphisms and of $E_0$-semigroups on
  factors}}
\end{center}

\bigskip
\begin{center}
Panchugopal Bikram\footnote{The Institute of Mathematical Sciences, Chennai}, 
Masaki Izumi\footnote{Graduate School of Science, Kyoto University, Kyoto},
R. Srinivasan\footnote{Chennai Mathematics Institute, Chennai},\\ and V.S. Sunder$^1$
\end{center}

\begin{abstract}
We examine what it means to say that certain endomorphisms  of a
factor (which we call equi-modular) are extendable. We obtain several
conditions for an equi-modular endomorphism to be extendable, and
single out one of them 
with a purely `subfactor flavour' as a theorem. We then exhibit the
obvious example of endomorphisms satifying the condition in this
theorem. We use our theorem to determine when every endomorphism in an
$E_0$-semigroup on a factor is extendable - which property is easily
seen to be a cocycle-conjugacy invariant of the $E_0$-semigroup. We 
conclude by giving examples of extendable $E_0-$semigroups, and by
showing that neither the Clifford flow on the hyperfinite $II_1$
factor, nor the free flow is extendable.
\end{abstract} 


 \bigskip






 



\section{Introduction}

We begin this note with a von Neumann algebraic version of the
elementary but extremely useful fact about being able 
to extend inner-product preserving maps from a total set of the domain
Hilbert space to an isometry defined on the entire domain. This leads
us to the notion of when a well-behaved (equi-modular, as we term it) 
endomorphism of a factorial probability space $(M,\phi)$ admits a
natural extension to an endomorphism of $L^2(M,\phi)$. After deriving
some equivalent conditions under which an endomorphism is extendable,
we exhibit examples of such extendable endomorphisms. 

We then pass to $E_0$-semigroups $\alpha = \{\alpha_t: t \geq 0\}$ of
factors, and observe that extendability of this semigroup (i.e.,
extendability of each $\alpha_t$) is a cocycle-conjugacy invariant of
the semigroup. We identify a necessary condition for extendability of
such an $E_0$-semigroup, which we then use to show that the Clifford
flow on the hyperfinite $II_1$ factor is not extendable. 

Our notion of extendable $E_0$-semigroups is
related to a notion called `regular semigroups' in \cite{ABS}, where
they erroneously claim to prove that the Clifford flow is extendable.

We start by setting up some notation. For any index set $I$, we write
$I^* = \bigcup_{n=0}^\infty I^n$ where $I^0 
= \emptyset$, and   $\i \vee \j = (i_1, \cdots ,i_m) \vee (j_1, \cdots
, j_n) =  (i_1, \cdots ,i_m, j_1, \cdots , j_n) $
whenever $\i = (i_1, \cdots ,i_m), \j = (j_1, \cdots , j_n) \in I^*$.

By a von Neumann probability space, we shall mean a pair $(M,\phi)$
consisting of a  von Neumann algebra and a normal state. 
For such an $(M,\phi)$, and an $x \in M$, we
shall write  $\hat{x} = \lambda_M(x) \widehat{1_M}$ and
$\widehat{1}_M$ for the cyclic  vector for $\lambda_M(M)$ in $L^2(M,\phi)$. 

Recall that the central support of the normal state $\phi$ is the central
projection $z (=:z_\phi)$ such that $ker(\lambda_M) = M(1-z)$. Clearly
$z_\phi = 1_M$ if $M$ is a factor.

Finally, if $\{x_i:i \in I\} \subset M$, and $\i =  (i_1, \cdots, i_n)
\in I^n$, we shall write $x_\i = x_{i_1} x_{i_2} \cdots x_{i_n}$. We also use $[S]$ either to denote 
the norm (respectively strong) closure of the span, for $S \subseteq \CH$ (respectively $S\subseteq \CL(\CH)$), for
any Hilbert space $\CH$.

\section{An existence result}
\begin{Proposition}\label{P}
Let $(M_i,\phi_i), i=1,2$ be von Neumann probability
spaces with $z_i = z_{\phi_i}$. Suppose $S^{(j)} = \{x^{(j)}_i:i \in I\}$ is a set of self-adjoint  
elements which generates $M_j$ as a von Neumann algebra, for
$j=1,2$. (Note the crucial assumption that both the $S^{(j)}$ are
indexed by the same set.)  Suppose 
\be\label{trcond}
\phi_1(x^{(1)}_\i)= \phi_2(x^{(2)}_\i) ~\forall \i \in I^*~.
\ee

Then there exists a unique isomorphism $\theta:M_1z_1 \rar M_2z_2$ such that
$\phi_2 \circ \theta|_{M_1z_1} = \phi_1|_{M_1z_1}$ and
$\theta(x^{(1)}_iz_1) = x^{(2)}_iz_2 ~\forall i \in I$. 
\end{Proposition}

\begin{proof}
The hypothesis implies that, for $j=1,2$, the set  $\{x^{(j)}_\i:\i
\in I^*\}$ linearly spans a *-subalgebra which is
necessarily $\sigma$-weakly dense in $M_j$.
Since $\langle \widehat{x^{(1)}_\i}, \widehat{x^{(1)}_\j} \rangle =  \langle \widehat{x^{(2)}_\i}, \widehat{x^{(2)}_\j} \rangle
~\forall \i,\j \in I^*$, there exists a unique unitary
operator $u:L^2(M_1,\phi_1) \rar L^2(M_2,\phi_2)$ such that $u
\widehat{x^{(1)}_\i} = \widehat{x^{(2)}_\i} ~\forall \i \in I^*$ . 

Now observe that
\beast
u\lambda_{M_1}(x^{(1)}_\i)u^*\widehat{x^{(2)}_\j} &=&
u\lambda_{M_1}(x^{(1)}_\i) \widehat{x^{(1)}_\j}\\
&=& u \widehat{x^{(1)}_{\i \vee \j} }\\
&=& \widehat{x^{(2)}_{\i \vee \j} }\\
&=& \lambda_{M_2}(x^{(2)}_\i)\widehat{x^{(2)}_{\j} }~;
\eeast
and hence that $u\lambda_{M_1}(x^{(1)}_\i)u^* =
\lambda_{M_2}(x^{(2)}_\i) ~\forall \i \in I$.

On the other hand, $\{x \in M_1: u\lambda_{M_1}(x)u^* \in
\lambda_{M_2}(M_2)\}$ is 
clearly a von Neumann subalgebra of $M_1$; since this has been shown to
contain each $x_\i^{(1)}$, we may deduce that this must be all of $M_1$. 
Now notice that $L^2(M_j,\phi_j) = L^2(M_jz_j,\phi_j|_{M_jz_j})$, that
$\lambda_{M_j}(x) = \lambda_{M_jz_j}(xz_j) ~\forall x \in M_j$, and that
$\lambda_{M_jz_j}$ maps $M_jz_j$ isomorphically onto its image.

The proof is completed by defining
\[\theta(x) = \lambda_{M_2 z_2}^{-1}(u\lambda_{M_1}(x)u^*) ~\forall x
\in M_1z_1.\]
\end{proof}

\begin{Remark}\label{end}
\ben
\item  In the proposition, even if it is the case that $N :=
  \{x^{(2)}_i:i \in I\}^{\prime\prime} \subsetneq M_2$,  we can still apply the
 result to $(N,\phi_2|_{N})$ in place of $(M_2,\phi_2)$ and deduce
 the existence of a normal homomorphism of 
$M_1$ into $M_2$ which sends $x^{(1)}_i$ to $x^{(2)}_iz$ for
each $i$ (and $1_{M_1}$ to the projection $z = z_{\phi_2|_N} \in N$).
\item In the special case that the $N$ of the last paragraph is a factor, the $z$ there is nothing but $id_{M_2}$ and in particular, Proposition \ref{P} can be strengthened as follows:

{\em Let $(M_j,\phi_j), j=1,2$ be von Neumann probability
spaces. Suppose $S^{(j)} = \{x^{(j)}_i:i \in I\} \subset M_j$ is a set of self-adjoint  elements such that $S^{(1)\prime \prime} = M_1$ 
and $S^{(2)\prime \prime}$ is a factor $N \subset M_2$. Suppose 
\be\label{trcond1}
\phi_1(x^{(1)}_\i)= \phi_2(x^{(2)}_\i) ~\forall \i \in I^*~.
\ee

Then there exists a unique normal $*$-homomorphism $\theta:M_1 \rar N \subset M_2$
such that $\theta(x^{(1)}_i) = x^{(2)}_i$ for all $i \in I$.} 
\een
\end{Remark}

\begin{Corollary}\label{faccor}
\ben
\item
If $\theta_i$ is a $\phi_i$-preserving unital endomorphism of a von
Neumann probability 
space $(M_i,\phi_i)$, for $i \in \Lambda$, then there exists:
\ben
\item a unique unital endomorphism $\otimes_{i \in \Lambda} 
\theta_i$ of the tensor product $(\otimes_{i \in \Lambda} M_i,
\otimes_{i \in \Lambda}\phi_i)$ such 
that $$(\otimes_{i \in \Lambda} \theta_i)(\otimes_{i \in \Lambda}
x_i) = z (\otimes_{i \in \Lambda} \theta_i(x_i)) ~\forall x_i = x_i^*
\in M_i;$$ 
\item a unique unital endomorphism $\ast_{i \in \Lambda} \theta_i$ of
the free product $(\ast_{i \in \Lambda} M_i, \ast_{i \in
  \Lambda}\phi_i)$ such that $$(\ast_{i \in
  \Lambda} \theta_i)(\lambda(x_j)) = z \lambda(\theta_j(x_j)) ~\forall
x_j \in M_j$$ where we simply write $\lambda$ for each `left-creation
representation' $\lambda:M_j \rar \CL( \ast_{i \in \Lambda} L^2(M_i,\phi_i))$ for
every $j \in I$.
\een
In the above existence assertions, the symbol  $z$ represents 
an appropriate projection (= image of the identity of the domain of the
endomorphism in question).

\item
If each $M_i$ above is a factor, then (the $z$ in the above statement
can be ignored, as it is the identity of the appropriate algebra) and
all endomorphisms above are unital monomorphisms.
\een 
\end{Corollary}

\begin{proof}
It is not hard to see that Remark \ref{end}(1) is applicable to $S^{(1)} =
\{\otimes_i x_i :  x_i = x_i^* \in M_i, x_i = 1_{M_i} \mbox{ for all
  but finitely many }i\}$  and $S^{(2)} = \{\otimes_i \theta_i(x_i) :
x_i = x_i^* \in M_i, x_i = 1_{M_i} \mbox{ for all but finitely many
}i\}$ (resp., $S^{(1)} = \{ \lambda(x_i): i \in \Lambda, x_i = x_i^*
\in M_i,   \phi_i(x_i) = 0\}$ and $S^{(2)} = \{ \lambda(\theta_i(x_i)) 
:  i \in \Lambda, x_i = x_i^* \in M_i, \phi_i(x_i)=0\}$.

The second fact follows from Remark \ref{end}(2) because normal
endomorphisms of factors are unital isomorphisms onto their images,
and the tensor (resp., free) product of factors is a factor.
\end{proof}

For later reference, the next lemma identifies the
central support $z_\phi$ of a normal state $\phi$ on a von Neumann
algebra in the simple special case when $\phi$ is a vector-state. 

\begin{Lemma}\label{vectst}
Suppose $N \subset \CL(\CH)$ is a von Neumann algebra, $\xi \in \CH$
is a unit vector, and $\phi$ is the vector state defined on $N$ by
$\phi(x) = \langle x\xi, \xi \rangle$. If $\CH_0 = \overline{N\xi}$,
then a candidate for `the GNS triple for $(N,\phi)$' is given by 
$(\CH_0, id_N|_{\CH_0},\xi)$. In particular, the central support of
$\phi$ is given by the projection $z = \wedge \{p \in N: ran~p \supset
\CH_0\}$ and $ran ~z = [N^\prime N \xi]$.
\end{Lemma}

\begin{proof}
It is clear that $\xi$ is a cyclic vector for $N|_{\CH_0}$ and the
assertion regarding GNS triples follows. 
Hence if $z \in \CP(Z(N))$ is such that $N(1-z) = ker~id_N|_{\CH_0}$,
then $z = \wedge\{p\in\CP(N): p|_{\CH_0} = (1_N)|_{\CH_0}\}$, i.e.,
$z = \wedge \{p \in \CP(N): ran~p \supset \CH_0\}$. As $z$ is the
smallest projection in $(N \cap N^\prime)$ whose range contains
$N\xi$, or equivalently the smallest subspace containing $[N\xi]$
which is invariant under $(N \cap N^\prime)^\prime$, equivalently
invariant under $N^\prime N$, the last assertion follows.
\end{proof}

\section{Extendable endomorphisms}

For the remainder of this paper, we make the standing assumption that
$\phi$ is a faithful normal state on a factor $M$. We identify $x \in M$ 
with $\lambda_M(x)$, and simply write $J$ and $\Delta$ for the modular conjugation
operator $J_\phi$ and the modular operator $\Delta_\phi$ respectively. 
Recall, thanks to the Tomita-Takesaki theorem that
$j = J(\cdot)J$ is a *-preserving conjugate-linear 
isomorphism of $\CL(L^2(M,\phi))$ onto itself, which maps $M$ and
$M^\prime$ onto one another,
and that $\widehat{1_M}$ is also a cyclic and separating vector for
$M^\prime$. We shall assume that $\theta$ is a normal unital
*-endomorphism which preserves $\phi$. The invariance assumption  
$\phi \circ \theta = \phi$ implies that there exists a unique isometry
$u_\theta$ on $L^2(M,\phi)$ such that $u_\theta x \widehat{1_M} = \theta(x)
\widehat{1_M}$ and equivalently, that $u_\theta x = \theta(x) u_\theta
~\forall x \in M$ and $u_\theta\widehat{1_M}=\widehat{1_M}$.

\begin{Definition}\label{equi-modular}
If $M,\phi, \theta$ are as above, and if the associated isometry
$u_\theta$ of $L^2(M,\phi)$ commutes with  
the modular conjugation operator $J (=J_\phi)$, we shall simply say
$\theta$ is a {\bf equi-modular} (as this is related to endomorphisms
commuting with the modular automorphism group) endomorphism of the  
factorial non-commutative probability space $(M, \phi)$. 
\end{Definition}

\begin{Remark}\label{iz1}
 It is true that if $\theta$ is an equi-modular
  endomorphism of a factor $M$ as above, then there always exists a
  $\phi$-preserving faithful normal conditional expectation $E:M \rar
  \theta(M)$, and in fact $u_\theta u_\theta^*$ is the Jones
  projection associated to this conditional expectation. 
  For this, notice to start with, that as $\theta$ is a *-homomorphism, $u_\theta$ commutes with the conjugate-linear 
  Tomita operator $S$ (which has $M\widehat{1_M}$ as a core and maps $x\widehat{1_M}$ to $x^*\widehat{1_M}$ for $x \in M$). 
  More precisely, we have $Su_\theta\supset u_\theta S$, meaning that whenever $\xi$ is in the domain of $S$, so is 
  $u_\theta\xi$ and $Su_\theta\xi=u_\theta S\xi$ holds. 
  Since $u_\theta$ commutes with $J$ and $S=J\Delta^{1/2}$, we have $\Delta^{1/2}u_\theta\supset u_\theta \Delta^{1/2}$, 
  and so $\Delta^{it}$ commutes with $u_\theta$ for any $t\in \R$. 
 Hence, conclude that for $x \in M$ and $t \in \R$, we have
 \beast
 \theta(\sigma^\phi_t(x)) \widehat{1_M} &=& u_\theta \Delta^{it} x \Delta^{-it} \widehat{1_M}
 = u_\theta \Delta^{it} x  \widehat{1_M}\\
 &=& \Delta^{it} u_\theta x \widehat{1_M}= \Delta^{it} \theta(x) \widehat{1_M}\\
 &=& \sigma^\phi_t(\theta(x)) \widehat{1_M}~.
 \eeast
 As $\widehat{1_M}$ is a separating vector for $M$, deduce that
 \[\theta(\sigma^\phi_t(x)) = \sigma^\phi_t(\theta(x)) ~\forall x \in
 M, t \in \R~.\] Hence $\sigma^\phi_t(\theta(M)) = \theta(M)
 ~\forall t\in \R$ and it follows from Takesaki's theorem (see \cite[Section 4]{Ta}) that there exists a unique $\phi$-preserving
 conditional expectation $E$ of $M$ onto the subfactor $P = \theta(M)$. 
 It is true, as the definition shows, that $e_\theta =
 u_\theta(u_\theta)^*$ is the orthogonal projection onto $[P\widehat{1_M}]$
 and $E(x)e_\theta = e_\theta x e_\theta ~\forall x \in M$.
\end{Remark}

\bt \label{prime-2}
Suppose $\theta$ is an equi-modular endomorphism of a factorial non-commutative probability space $(M, \phi)$.  Then,

\ben
\item The equation $\theta^\prime = j \circ \theta \circ j$
  defines a unital normal *-endomorphism of $M^\prime$ which preserves
$\phi^\prime = \overline{\phi \circ j}$; and
\item We have an identification
\beast L^2(M',\phi ') &=& L^2(M,\phi)\\
\widehat{1_{M'}} &=& \widehat{1_M}\\
u_{\theta '} & =& u_\theta
\eeast
\item there exists a unique endomorphism
  $\theta^{(2)}$ of $\CL(L^2(M,\phi))$ satisfying
\[\theta^{(2)}(xj(y)) = \theta(x) j(\theta(y))z ,
~\forall x,y \in M\]
where $z = \wedge \{p \in (\theta(M) \cup
\theta^\prime(M^\prime))^{\prime\prime}: ran(p) \supset
\{\widehat{\theta(x)}:x \in M\}\}$.
\een
\et

\begin{proof}
\ben
\item
It is clear that $\theta^\prime = j \circ \theta \circ j$ is a unital
normal linear *-endomorphism of $M^\prime$ and that 
\[\overline{\phi^\prime \circ \theta^\prime} = \overline{\phi^\prime} \circ \theta^\prime = (\phi \circ j) \circ (j \circ
\theta \circ j) = (\phi \circ \theta) \circ j = \phi \circ j =
\overline{\phi^\prime}~,\] 
thereby proving (1).

\item This follows from the facts that $\widehat{1_M}$ is a cyclic and separating vector for $M$ and hence also for $M'$, the definition of $\phi '$ which guarantees that
\beast
\langle j(x) \widehat{1_{M'}}, j(y) \widehat{1_{M'}} \rangle &=& \phi '(j(y)^*j(x))\\
&=& \phi '(j(y^*x))\\
&=& \overline{\phi(y*x)}\\
&=& \phi(x^*y)\\
&=& \langle y \widehat{1_M}, x \widehat{1_M} \rangle\\
&=& \langle Jx \widehat{1_M}, Jy \widehat{1_M} \rangle\\
&=& \langle JxJ \widehat{1_M}, JyJ \widehat{1_M} \rangle\\
&=& \langle j(x) \widehat{1_M}, j(y) \widehat{1_M} \rangle
\eeast
and the definitions of the `implementing isometries', which show that
\beast
u_{\theta '}(j(x)\widehat{1_{M'}}) &=& \theta '(j(x)) \widehat{1_{M'}}\\
&=& j(\theta(x)) \widehat{1_{M'}}\\
&=& J \theta(x) J \widehat{1_{M}}\\
&=& J \theta(x) \widehat{1_{M}}\\
&=& J u_\theta x \widehat{1_{M}}\\
&=& u_\theta Jx \widehat{1_{M}}\\
&=& u_\theta JxJ \widehat{1_{M}}\\
&=& u_\theta j(x) \widehat{1_{M'}}~.
\eeast

\item Notice that if $x,y \in M$, then
\beast
\langle \theta(x) J \theta(y)J  \widehat{1_M} ,
\widehat{1_M}\rangle 
&=& \langle \theta(x) J \theta(y)  \widehat{1_M} ,
\widehat{1_M}\rangle \\
&=& \langle \theta(x) J u_\theta y \widehat{1_M} ,
\widehat{1_M}\rangle \\
&=& \langle \theta(x) u_\theta J y \widehat{1_M} ,
\widehat{1_M}\rangle \\
&=& \langle u_\theta x J y \widehat{1_M} ,
\widehat{1_M}\rangle\\
&=& \langle u_\theta x J y \widehat{1_M} ,
u_\theta \widehat{1_M}\rangle\\
&=& \langle x J y J \widehat{1_M} ,
\widehat{1_M}\rangle ~,
\eeast
where we have used the fact that $\theta$ is equi-modular.

Set $S^1 = \{xj(y): x = x^*, y=y^*, x,y \in M\}$, and
$S^{(2)} = \{\theta(x)j(\theta(y)) : xj(y) \in S^{(1)}\}$, and deduce
from the factoriality of $M$
that  $S^{(1)\prime \prime} = \CL(L^2(M,\phi))$.

Now we wish to apply  Remark \ref{end}(1) with $N =
S^{(2)\prime\prime} = \theta(M) \vee j(\theta(M))$  (where, both here
and in the sequel, we write $A \vee B = (A \cup B)^{\prime \prime}$
for the von Neumann algebra generated by von Neumann algebras $A$ and
$B$) and $\phi_1 = \phi_2 = \langle (\cdot) \hat{1}_M, \hat{1}_M
\rangle$. For this, deduce  from Lemma \ref{vectst} that 
\beast
z &=& \wedge \{p \in \CP(N): ran~p \supset N\hat{1_M}\}\\
&=& \wedge \{p \in \CP(N): ran~p \supset \{\theta(x)\hat{1_M},  \theta^\prime(j(x)) \hat{1_M}: x \in M\}\}\\
&=& \wedge \{p \in \CP(N): ran~p \supset \{\widehat{\theta(x)}:x \in M\}\}\\
\eeast
and the proof of the Theorem is complete.
\een
\end{proof}

\begin{Remark}\label{iz2}
It must be observed that the projection $z$ of Theorem \ref{prime-2} is nothing but the central support of
 the projection $e_\theta = u_\theta u_\theta^*$ in $P^\prime \cap P_1$
 where $P = \theta(M) \subset M \subset P_1$ is Jones' basic
 construction (thus, $P_1 = JP^\prime J$) since, by Lemma \ref{vectst}, we have: 
\[ran ~z = [(P\vee JPJ)'(P\vee JPJ)\widehat{1_M}]= [(P'\cap P_1)e_\theta L^2(M,\phi)]~.\] 
This is because
\beast
[(P\vee JPJ)\widehat{1_M}] &=&[PJPJ\widehat{1_M}]=[PJP\widehat{1_M}]\\
&=& [PJu_\theta M \widehat{1_M}]= [Pu_\theta JM \widehat{1_M}]\\
&=& [Pu_\theta M \widehat{1_M}]~=[P\widehat{1_M}]~.
\eeast 
In particular, since $e_\theta$ is a minimal projection in $P' \cap P_1$, 
its central support $z$ in $P'\cap P_1$ is 1 if and only if $P' \cap P_1$ is a type I factor. 
In the following corollary, we continue to use the symbols
 $P$ and $P_1$ with the  meaning attributed to them here. 
\end{Remark}

The following corollary is an immediate consequence of Lemma
\ref{vectst}, Theorem \ref{prime-2} and Remark \ref{iz2}.

\begin{Corollary}\label{extendable}
Let $\theta$ be a equi-modular endomorphism of a factorial
non-commutative probability space $(M, \phi)$ in standard
form (i.e., viewed as embedded in $\CL(L^2(M,\phi))$ as above).
The following conditions on $\theta$ are equivalent:
\ben 
\item there exists a unique unital normal $*$-endomorphism $\theta^{(2)}$ of
  $\CL(L^2(M,\phi))$ such that $\theta^{(2)}(x) = \theta(x)$ and
  $\theta^{(2)}(j(x)) = j(\theta(x))$ for all $x \in M$. 
\item $P \vee JPJ $ is a factor; and in this
  case, it is  necessarily a type $I$ factor.
\item $(P\vee JPJ)'=P'\cap P_1$ is a factor; and in this
  case, it is  necessarily a type $I$ factor.
\item $\{x\widehat{y}: x \in P'\cap P_1, y \in P\}$ is
  total in $L^2(M,\phi)$. 
\een

An endomorphism of a factor which satisfies the equivalent conditions
above will be said to be 
{\bf extendable}.
\end{Corollary}

\begin{Remark}\label{phiext}
It should be noted that extendability is not a property of just an
endomorphism $\theta$ but it is also dependent on a state which is not
only left invariant under the enomorphism but must also satisfy the
requirement we have called equi-modular. Strictly speaking, we should
probably talk of $\phi$-extendability, but shall not do so in the
interest of notational convenience.
\end{Remark}

\begin{Theorem}\label{rela} Let the notation be as above. 
Then the following conditions are equivalent: 
\ben
\item $\theta$ is extendable. 
\item $M=(M\cap \theta(M)')\vee \theta(M)$. 
(Note that the right-hand side is naturally identified with the von Neumann algebra 
tensor product $(M\cap \theta(M)')\otimes \theta(M)$ in this case.)
\een
\end{Theorem}

\begin{proof} Recall that $P=\theta(M)$ is globally preserved by the modular automorphism group 
$\{\sigma^\phi_t\}_{t\in \R}$ and there exists a $\phi$-preserving faithful normal conditional expectation 
$E$ from $M$ onto $P$. 
Thus $(M\cap P')\vee P$ is naturally identified with the von Neumann algebra 
tensor product $(M\cap P')\otimes P$ (see \cite[Corollary 1]{Ta}). 
If we assume the second condition in the statement, the basic construction for $P\subset M$ essentially 
comes from that of $\C\subset (M\cap P')$, and so $P'\cap P_1$ is a type I factor. 
This means that $\theta$ is extendable. 

Assume that $\theta$ is extendable now. 
We will show that $Q:=(M\cap P')\vee P$ coincides with $M$. 
For this, it suffices to show $[Q \widehat{1_M}]=L^2(M,\phi)$. 
Indeed, since $P$ is globally preserved by $\sigma^\phi_t$, so is $Q$, and there exists a $\phi$-preserving 
faithful normal conditional expectation from $M$ onto $Q$ thanks to Takesaki's theorem. 
Thus if $Q$ were a proper subalgebra of $M$, $[Q \widehat{1_M}]$ would be a proper subspace of 
$L^2(M,\phi)$. 

Let $\hat{E}$ be the dual operator valued weight from $P_1$ to $M$ (see \cite{Ko} for the definition of 
$\hat{E}$ and its properties). 
Since $E\circ \hat{E}(e_\theta)=1<\infty$ and $P'\cap P_1$ is a factor, the restriction of 
$E\circ \hat{E}$ to the type I factor $P'\cap P_1$ is a faithful normal semifinite weight (see \cite[Lemma 2.5]{ILP}). 
Thus there exists a (not necessarily bounded) non-singular positive operator $\rho$ affiliated to $P'\cap P_1$ 
satisfying $\sigma^{E\circ \hat{E}}_t=\mathrm{Ad} \rho^{it}$ and 
$$E\circ \hat{E}(a)=\lim_{n\to \infty}\mathrm{Tr}(\rho(1+\frac{1}{n}\rho)^{-1} a),\quad \forall 
a\in (P'\cap P_1)_+,$$
where $\{\sigma^{E\circ \hat{E}}_t\}_{t\in \R}$ is the relative modular automorphism group (the restriction 
of $\{\sigma^{\phi\circ \hat{E}}_t\}_{t\in \R}$ to $P'\cap P_1$). 
Note that the trace $\mathrm{Tr}$ makes sense as $P'\cap P_1$ is a type I factor. 
 
From the above argument we see that there exists a partition of unity $\{e_i\}_{i\in I}$ consisting of 
minimal projections $e_i\in P'\cap P_1$ with $E\circ \hat{E}(e_i)<\infty$. 
Since $e_\theta$ is a minimal projection in $P'\cap P_1$ satisfying 
$\sigma^{E\circ \hat{E}}_t(e_\theta)=e_\theta$ and $E\circ \hat{E}(e_\theta)=1$, 
we may assume $0\in I$ and $e_0=e_\theta$. 
Let $\{e_{ij}\}_{i,j\in I}$ be a system of matrix units in $P'\cap P_1$ satisfying $e_{ii}=e_i$. 
Then we can apply the push down lemma \cite[Proposition 2.2]{ILP} to $e_{0i}$, and we have 
$e_{0i}=e_\theta q_i$, where $q_i=\hat{E}(e_{0i})\in P'\cap M$. 
Now for any $x\in M$, we have 
$$x\widehat{1_M}=\sum_{i\in I}e_{ii}x\widehat{1_M}=\sum_{i\in I}q_i^*e_\theta q_i x\widehat{1_M}=\sum_{i\in I}q_i^*E(q_i x)\widehat{1_M},$$
which shows $[Q\widehat{1_M}]=L^2(M,\phi)$. 
\end{proof}

\section{Examples of Extendable Endomorphisms}
Note that any automorphism on a factor is extendable, since the conditions in Corollary \ref{extendable} are satisfied. 

Let $\CR$  denote  the hyperfinite $II_1$ factor and $M$ be any 
$II_1$ factor which is also a  McDuff factor; i.e.,  $M \otimes \CR \cong M$.
Let $\alpha : M \otimes \CR \mapsto M $ be an isomorphism and 
$\beta : M \mapsto  M \otimes \CR $ be the monomorphism defined by 
$ \beta ( m ) = m \otimes 1 $, for $ m\in M$.  Let us write 
$ \theta = \beta \circ \alpha$. so $\theta $ is an endomorphism 
of $M \otimes \CR$ such that 
$\theta (M \otimes \CR ) = M \otimes 1$. As $M \otimes \CR$
is a $II_1$ factor,  the endomorphism $\theta$ is necessarily
equi-modular. Now by corollary  
\ref{extendable}, showing  that $\theta$ is  extendable is equivalent to showing that $\{ \theta ( M \otimes \CR ) \vee J  \theta( M \otimes \CR ) J \} $ is a
type $I$ factor, where $J$ is the modular conjugation of $M
\otimes\CR$, which, of course, is $J_M \otimes J_R$. Note that
\beast
\{ \theta ( M \otimes \CR ) \vee J  \theta( M \otimes \CR ) J \} &=& 
\{ M\otimes 1 \vee J   (M \otimes 1)  J \}\\
&=& \{ M\otimes 1 \vee J_M M J_M\otimes 1 )  \}\\
&=& \CL(L^2(M) \otimes 1 
\eeast 
So $\{ \theta ( M \otimes \CR ) \vee J  \theta( M \otimes \CR ) J \} $ is 
a type $I$ factor. That is $\theta$ is extendable. 

\section{Extendability for $E_0$-semigroups}
\begin{Definition}
$\{\alpha_t: t \geq 0\}$ is said to be an $E_0$-semigroup on a von
Neumann probability space $(M,\phi)$ if:
\ben
\item $\alpha_t$ is a $\phi$-preserving normal unital *-homomorphism
  of $M$ for each $t \geq 0$;
\item $\alpha_0 = id_M$and $\alpha_s \circ \alpha_t = \alpha_{s+t}$; and
\item $[0,\infty) \ni t \mapsto \rho(\alpha_t(x))$ is continuous for
  each $x \in M, \rho \in M_*$.
\een

Suppose $\alpha_t$ is (equi-modular and) extendable
for each $t$, then we say that the $E_0$-semigroup
$\alpha$ is extendable.
\end{Definition}

\begin{Remark}\label{phiext1}
A remark along the lines of Remark \ref{phiext}, with endomorphism
replaced by $E_0$-semigroup, is in place here.
\end{Remark}

\begin{Proposition}\label{exte0}
Suppose $\alpha = \{\alpha_t: t \geq 0\}$ is an
$E_0$-semigroup on a factorial non-commutative
 probability space $(M,\phi)$.
\ben
\item The equation $\alpha_t^\prime(x^\prime) =
  j(\alpha_t(j(x^\prime))$ defines an $E_0$-semigroup on
  $(M^\prime,\phi^\prime)$,
where $\phi^\prime(x^\prime) = \omega_{\widehat{1_M}}(x^\prime) =
\langle x^\prime \widehat{1_M}, \widehat{1_M} \rangle$;
\item If $\alpha$ is extendable
for each $t$, then there exists a unique $E_0$-semigroup $\{\alpha^{(2)}_t: t \geq
  0\}$ on $(\CL(L^2(M,\phi)),\omega_{\widehat{1_M}})$ such that
  $\alpha^{(2)}_t(xx^\prime) = \alpha_t(x)\alpha_t^\prime(x^\prime)
  ~\forall x \in M, x^\prime \in M^\prime$.
\een
\end{Proposition}

\begin{proof} Existence of the endomorphisms $\alpha^\prime_t$ and
  $\alpha^{(2)}_t$ is guaranteed by Corollary \ref{extendable}, 
The equation $\alpha^\prime_t = j \circ \alpha_t \circ j$ shows that
$\{\alpha^\prime_t: t \geq 0\}$ inherits the property of being an
$E_0$-semigroup from that of $\{\alpha_t: t \geq 0\}$. The corresponding property for $\{\alpha^{(2)}_t: t \geq 0\}$ is now
seen to follow easily from the uniqueness assertion in Corollary
\ref{extendable}(1).
\end{proof}

By using standard arguments from the theory of $E_0$-semigroups on
type I factors, we can strengthen Corollary \ref{extendable} in the
case of $E_0$-semigroups thus:

\begin{Proposition}\label{extn}
Let $\alpha = \{\alpha_t: t \geq 0\}$ be an
$E_0$-semigroup on a factorial non-commutative 
probability space $(M,\phi)$, and suppose $\alpha_t$ is equi-modular 
for each $t$. Suppose $M$ is 
acting standardly on $ \CH = L ^{2} ( M )$. Consistent with the
notation of Remark \ref{iz2}, we shall write $P(t) = \alpha_t(M)
\subset M \subset P_1(t)$ for Jones' basic construction.

The following conditions on $\alpha$ are equivalent.
\ben
\item $\alpha$ is extendable.
\item $P'(t) \cap P_1(t)$ is either: (i) a factor of type $I_1$ (i.e., is
  isomorphic to $\C$) and $\alpha_t$ is an automorphism for all $t$;
  or (ii) a factor of type $I_\infty$ for all $t$ and  no $\alpha_t$
  is an automorphism. 
\een
\end{Proposition}
\begin{proof}
\bigskip
$(1) \Rar (2)$ If each $\alpha_t$ is extendable, then  $P'(t) \cap
P_1(t)$ is a factor of type $I_{n_t}$, say, by Corollary
\ref{extendable}. The first fact to be noted is that if
$\{\CE_{\alpha^{(2)}}(t): t \geq 0\}$ is the product system associated
to the $E_0$-semigroup $\alpha^{(2)}$ on $\CL(L^2(M))$, then $n_t$ is the
dimension of the Hilbert space $\CE_{\alpha^{(2)}}(t)$. Hence $\{n_t:t
\geq 0\}$ is a multiplicative semi-group of integers. Hence either
$n_t$ is constant in $t$ (identically 1 or identically infinity).

$(2) \Rar (1)$ follows from Corollary \ref{extendable}.
\end{proof}

\begin{Remark} Let $\alpha = \{\alpha_t: t \geq 0\}$ be an
$E_0$-semigroup on a factorial non-commutative 
probability space $(M,\phi)$, and suppose $\alpha_t$ is equi-modular 
for each $t$.
If  $\alpha_t$ is an extendable endomorphism for some $t>0$, then $\alpha_s$ is also extendable for all $0\leq s<t$.
Indeed, $\alpha_t$ being extendable means that $M$ as a $P(t)-P(t)$ bimodule is a direct  sum of copies of $P(t)$, and hence $P(t-s)$ as a $P(t)-P(t)$ bimodule is also a direct sum of copies of $P(t)$.  This means that $\alpha_{t-s}(M)$ is generated by $\alpha_{t-s}(M)\cap \alpha_t(M)^\prime$Êand $\alpha_t(M)$, which means  $\alpha_s$ is extendable. Now, since the compositions of extendable endomorphisms are extendable, the $E_0-$semigroup $\alpha$ itself is extendable.
\end{Remark}

Now let us consider the following spaces;
\[
 E^{\alpha_t} = \{ T \in \CL(L^2(M)) : \alpha_t(x)T = Tx, \text{ forall }x \in M \};
\]
\[
 E^{\alpha_t'} = \{ T \in \CL(L^2(M)) : \alpha_t'(x')T = Tx', \text{ forall }x' \in M' \}.
\]
For every $t \geq 0 $, we write $ H(t) =  E^{\alpha_t} \cap E^{\alpha_t'}$. Then 
we have the following Lemma. 
\begin{Proposition}\label{prodlem}
Let $\alpha = \{\alpha_t: t \geq 0\}$ be an
$E_0$-semigroup on a factorial non-commutative 
probability space $(M,\phi)$ and suppose $\alpha_t$ is equi-modular 
for each $t$. 
If  $\alpha$ is extendable then 
 $$H= \{ (t , T ) : t \in (0,\infty) , T \in H(t) \}$$ is a product 
system (in the sense of \cite{Arv}) with the family of unitary maps 
$u_{st} : H(t) \otimes H(s) \mapsto H({s+t})$,
given by $$u_{st} ( T \otimes S) = TS ~\forall T \in H(t), S \in H(s).$$   
\end{Proposition}

\begin{proof}
Let $\alpha^{(2)} = \{ \alpha_t^{(2)} : t > 0 \}$
be the extension of $\alpha $ on $\CL(L^2(M))$. 
For $t > 0$, consider
\[
\CE(t) = \{ T \in \CL(L^2(M)) : \alpha^{(2)}_t(x) T = Tx ,\text{ for all } x \in \CL(L^2(M)) \}.
\]
We shall write $\CE = \{ (t, T ) : T \in \CE(t) \}$; then
$\CE$ is a product system (see \cite{Arv}), and  $H(t)= \CE(t)$ for every $t>0$. Indeed,
if $ T \in H(t) = E^{\alpha_t} \cap E^{\alpha_t^\prime} $, then
$\alpha_t(m)T = Tm$ for all $m \in M$ and 
$\alpha_t^\prime(m^\prime )T = Tm^\prime $ for all $m^\prime \in M^\prime $.
So it is clear that $\alpha_t^{(2)}(x)T = Tx $ for all $x \in M \cup
M'$ and hence also for all
$x \in (M \vee M') = \CL(L^2(M))$.
So, $T \in \CE(t)$, and  $H(t) \subset \CE(t)$. The reverse inclusion is
immediate from the definition $\alpha_t^{(2)}$. So we have 
$H(t) = \CE(t)$ and clearly $H$ is a product system.
\end{proof}

Now recall that an $E_0$-semigroup $\{\beta_t:
t \geq 0\}$ of a von Neumann probability space $(M,\phi)$ is said to
be a {\bf cocycle perturbation} of an $E_0$-semigroup $\{\alpha_t: t
\geq 0\}$ if there exists a weakly continuous family 
$\{u_t: t \geq 0\}$ of unitary elements of $M$ such that
\ben
\item $u_{t+s} = u_s \alpha_s(u_t)$; and
\item $\beta_t(x) = u_t \alpha_t(x) u_t^*$ for all $x \in M$ and $s,t
  \geq 0$.
\een
In such a case, we shall simply write
\[\{u_t: t \geq 0\}: \{\alpha_t : t \geq 0\} \simeq \{\beta_t: t \geq
0\}. \]

\begin{Proposition}\label{ccprop}
Suppose $\beta = \{\beta_t : t \geq 0\}$ is an  $E_0$
semigroup on a factorial probability space $(M,\phi)$, which is a
cocycle perturbation of another $E_0$
semigroup $\alpha = \{\alpha_t : t \geq 0\}$ on $(M,\phi)$ with
$\{u_t: t \geq 0\}: \{\alpha_t : t \geq 0\} \simeq \{\beta_t: t \geq
0\}. $
Then 
\ben
\item $\{j(u_t): t \geq 0\}: \{\alpha^\prime_t : t \geq 0\} \simeq
\{\beta^\prime_t: t \geq 0\}.$
\item If each $\alpha_t$ is extendable, as is each $\beta_t$, then
\[\{u_tj(u_t): t \geq 0\}: \{\alpha^{(2)}_t : t \geq 0\} \simeq
\{\beta^{(2)}_t: t \geq 0\}. \]
\een
\end{Proposition}

\begin{proof} The verifications are elementary and a routine
  computation. For example, once (1) has been verified, the
  verification of (2) involves such straightforward computations as:
if we let $U_t = u_t u^\prime_t$, where we write $u_t' = j(u_t)$, and
if $x \in M, x^\prime  \in M^\prime$, then
\beast U_t\alpha_t(x)U_t^* &=& u_t \alpha_t(x)  \u_t^* = \beta_t(x)\\
U_t\alpha^\prime_t(x^\prime)U_t^*  &=& u^\prime_t \alpha^\prime_t(x^\prime)u^{\prime *}_t = \beta^\prime_t(x^\prime)\\
\beta_t(M) \vee \beta^\prime_t(M^\prime) &=& U_t(\alpha_t(M) \vee \alpha^\prime_t(M^\prime))U_t^*\\
U_{s+t} &=& u_{s+t} u^\prime_{s+t}\\
&=& u_s \alpha_s(u_t) u^\prime_s \alpha^\prime_s(u^\prime_t)\\
&=& U_s \alpha^{(2)}_s(U_t)\\
\eeast
and
\beast
\beta^{(2)}_t(xx^\prime) &=& \beta_t(x) \beta^\prime_t(x^\prime)\\
&=& u_t \alpha_t(x) u_t^* u^\prime_t \alpha^\prime_t(x^\prime)
u^{\prime *}_t\\
&=& U_t \alpha^{(2)}_t(xx^\prime) U_t^*.
\eeast
\end{proof}

Recall that two $E_0$-semigroups $\{\alpha_t: t \geq 0\}$ and
$\{\beta_t: t \geq 0\}$ on a von Neumann algebra $M$ are said to be
{\bf conjugate} if there exists an automorphism $\theta$ of $M$ such
that $\beta_t \circ \theta = \theta \circ \alpha_t ~\forall t$, while
they are said to be {\bf cocycle conjugate} if each is conjugate to a
cocycle perturbation of the otheâçr.

\begin{Remark}\label{inddefandcc}
While the index of $E_0$-semigroups of type $I_\infty$ factors has
been well-defined, we may now define the index of an extendable $E_0$ 
semigroup $\alpha$ of an arbitrary factor as the index of
$\alpha^{(2)}$;
and we may infer from Proposition \ref{ccprop} that the index of an
extendable
$E_0$-semigroup of an arbitrary factor is invariant under cocycle
conjugacy - in the restricted sense that cocycle conjugate 
extendable $E_0$-semigroups have the same index. (One has to exercise
some caution here in that there is a problem with invariance of
equimodularity  under
cocycle conjugacy!)\footnote{We wish to thank the referee for pointing
  this out, which also led to the insertion of the Remarks
  \ref{phiext} and \ref{phiext1}.}
It is to be noted from Corollary \ref{extendable} that the extendability of an
$E_0$-semigroup, each of whose endomorphisms is equi-modular, is a
property which is invariant under cocycle conjugacy within 
the class of such $E_0$-semigroups. 
\end{Remark}

\begin{Proposition}\label{tpext}
If $\alpha = \{\alpha_t:t \geq 0\}$ (resp., $\beta = \{\beta_t:t \geq
0\}$) is an extendable $E_0$-semigroup of a factor $M$ (resp., $N$),
then $\alpha \otimes \beta = \{\alpha_t \otimes \beta_t:t \geq 0\}$ is an
extendable $E_0$-semigroup of the factor $M \otimes N$, and in fact,
\[(\alpha \otimes \beta)^{(2)} = \alpha^{(2)} \otimes \beta^{(2)}.\]
\end{Proposition}

\begin{proof}
The hypothesis is that $\alpha_t(M) \vee J \alpha_t(M) J$ and
$\beta_t(N) \vee J_N \beta_t(N) J_N$ are factors, for each $t \geq 0$,
while the conclusions follow from the definition of $\alpha \otimes \beta$.
\end{proof}

\section{Examples} 
First we give examples of extendable $E_0-$semigroups. Throughout this
section, let $\CH = L^2 (0, \infty ) \otimes \CK$, be the real Hilbert
space of square integrable functions taking values in is a real
Hilbert space $\CK$. We always denote by 
$(\cdot)_\C$ the complexification of $(\cdot)$. 
Let $S_t$ be the shift semigroup on $\CH_\C$ defined by
\begin{eqnarray*}(S_tf)(s) & = & 0, \quad s<t,\\
& = & f(s-t), \quad s \geq t.
\end{eqnarray*} Thus $(S_t:t \geq 0)$ is a  semigroup of isometries,
and we denote its restriction to $\CH$ also by $\{S_t \}$.

For the first set of examples, given by `canonical commutation
relations', we only need complex Hilbert spaces. Let $A\geq 1$ be a
complex linear operator on $\CH_\C$ such that $T =\frac{1}{2}(A-1)$ is
injective. Consider the the quasi free state on the CCR algebra over
$\CH_\C$ given by $$\varphi_A(W(f))=e^{-\frac{1}{2}\langle Af,
  f\rangle}= e^{-\frac{1}{2}\|\sqrt{1+2T}f\|^2} ~~\forall ~~f \in
\CH_\C.$$ The space underlying the corresponding GNS representation
may be identified with $\Gamma_s(\CH_\C)\otimes \Gamma_s(\CH_\C)$,
with the GNS representation being described by $$\pi(W(f) ) = W_0
(\sqrt{1+T}f)\otimes W_0(j\sqrt{T}f)~~\forall ~~f \in \CH_\C,$$ where
$\Gamma_s(\cdot)$ is the symmetric Fock space, $W_0(\cdot)$ is the
Weyl operator on $\Gamma_s(\CH_\C)$ and $j$ is an anti-unitary on
$\CH_\C$ induced by an anti-unitary operator on $\CK_\C$. The vacuum
vector $\Omega \otimes \Omega \in \Gamma_s(\CH_\C) \otimes
\Gamma_s(\CH_\C)$ is the cyclic and separating vector for
$M_A=\left\{\pi (W(f)): f \in \CH_\C\right\}^{\prime \prime}$ (see \cite{Ark}).

\begin{Example} 
Let $A=\frac{1+\lambda}{1-\lambda}$ with $\lambda\in (0,1)$, then
it is well-known that $M_\lambda=M_A$ is a type
$III_\lambda$ factor. 
There exists a unique $E_0-$semigroup $\beta_t^\lambda$ on $M_\lambda$ satisfying $$\beta^\lambda_t(\pi(W(f)))=\pi(W(S_tf))~~ \forall ~ f \in \CH_\C.$$ 

Further, $\{\beta_t^\gamma;t\geq 0\}$ is equi-modular and 
the relative commutant is given by $$\beta_t^\lambda(M_\lambda)^\prime\cap M_\lambda= \left\{\pi (W(f)): f \in (L^2 (0, t) \otimes \CK\right)_\C\}^{\prime \prime}.$$ Now theorem  \ref{rela} imply that  all these $E_0-$semigroups on type $III_\lambda$ factors are extendable. (See \cite{OV}, where these examples are discussed in more detail.)
\end{Example}

We will write $\CF(\CH)$ for the anisymmetric Fock space; thus
\[
\CF(\CH_\C) = \C\Omega \oplus  \CH_\C \oplus
\left(\CH_{\C}\wedge\CH_{\C} \right)\oplus
\left(\CH_{\C}\wedge\CH_{\C}\wedge\CH_{\C}\right) \oplus \cdots,
\]
where $\Omega$ is a fixed complex number with modulus 1.


Recall the left creation operator on $\CF(\CH)$ (corresponding to $f
\in \CH_{\C}$ given by
\beast a(f) \Omega &=& f\\ 
 a(f) ( \xi_{1} \wedge \xi_{2}\wedge \cdots \wedge \xi_{n}) &=& f\wedge\xi_{1} \wedge \xi_{2}\wedge \cdots \wedge \xi_{n} \text{ , } 
\xi_{i}\in \CH_{\C}~.
\eeast
 
These operators obey the {\em canonical anticommutation relations}:
\[
 a(f)a(g) + a(g)a(f) = 0 \text{ , } a(f)a(g)^{*} + a(g)^{*}a(f) =  \langle f\text{ , }g \rangle id_{\CF(\CH)}
\]
for all $ f \text{ , }  g \in \CH_{\C}$ , 

For any $f \in \CH$, let $u(f) = a(f) + a(f)^{*}$.
It is well-known that
the von Neumann  algebra $$\{ u(f) : f \in \CH \}^{''}\subseteq 
\CL(\CF(\CH_{\C}) )$$  
is the hyperfinite $II_{1}$ with cyclic and separating ( trace)vector $\Omega$. 

\begin{Example}
For every $t\geq 0$ there exist 
 a unique normal, unital  $*$-endomorphism $ \alpha_t : \CR \mapsto
 \CR $ satisfying 
\[
\alpha_t ( u(f) ) =  u(S_tf)) ~\forall f \in \CH_\C~.
\]
(Although this is a well-known fact, we remark that this is in
fact a consequence of Remark \ref{end} (2).)
Then $ \alpha = \{ \alpha_t : t \geq 0 \}$ is an 
$E_0$-semigroup on $\CR$, called the {\bf Clifford flow of
 rank dim $\CK$}.

It is known from \cite{Alev} that 
\be\label{cliffrc} \alpha_t(M)' \cap M = \{u(f)u(g): spt(f), spt(g)\subset [0,t]\}''~. \ee
 It follows from equation \ref{cliffrc} that if $spt(f) \subset [0,t]$, then $u(f) \Omega \perp \{(\alpha_t(\CR)' \cap \CR)\Omega \cup \alpha_t(\CR) \Omega\}$; (in fact the same assertion holds for any $a(f_1) \cdots a(f_{2n+1} \Omega$ for any $n$ and any $f_1,\cdots, f_{2n+1}$ with support in $[0,t]$.) Consequently, in view of $\Omega$ being a separating vector for $\CR$,  it is an easy  consequence of Theorem \ref{rela} that the
{\bf Clifford flow on $\CR$ (of any rank) is not extendable}.
\end{Example}

The Clifford flows of the hyperfinite II$_1$ factor are closely
related to another family  of $E_0-$semigroups, called the CAR
flows. (We should remember these are CAR flows on type II$_1$ factors,
not to be confused with the usual CAR flows on the type $I$ factor of all bounded operators on the antisymmetric Fock space.)  We recall the definition of  CAR algebra and some facts regarding the GNS representations of CAR algebras given by quasi-free states. 

For a complex Hilbert space $K$, the associated CAR algebra $CAR(K)$ is the universal $C^*-$algebra generated by a unit $1$ and elements $\{b(f):~f\in K\}$, subject to the following relations
\begin{itemize}
\item[(i)] $b(\lambda f)= \lambda b(f)$,
\item[(ii)] $b(f)b(g)+b(g)b(f)=0$,
\item[(iii)] $b(f)b^*(g)+b^*(g)b(f) =\langle f,g\rangle 1$,
\end{itemize}
for all $\lambda\in\C$, $f,g\in K$, where $b^*(f)=b(f)^*$.

Given a positive contraction $A$ on $K$, there exists a unique quasi-free state $\omega_A$ on $CAR(K)$ satisfying  $$\omega_A(b(x_n) \cdots b(x_1)b(y_1)^* \cdots b(y_m)^* ) = \delta_{n,m}det (\langle Ax_i, y_j\rangle),$$ where $det(\cdot)$ denotes the determinant of a matrix. Let $(H_A, \pi_A, \Omega_A)$ be the corresponding GNS triple. Then  $M_A=\pi_A(CAR(K))^{\prime\prime}$ is a factor.  

Here onwards we fix the contraction with $A = \frac{1}{2}$, then $M_A=\CR$ is the hyperfinite type II$_1$ factor and $\omega_A$ is a tracial state. We define the CAR flow on $\CR$ as follows.  

Now let $K=\CH_\C$. Then there exists a unique $E_0-$semigroup $\{\alpha_t\}$ on $\CR$ satisfying $$\alpha_t(\pi(b(f)))=\pi(b(S_tf))~~ \forall ~f \in \CH_\C.$$ This $\alpha$ is called  as the CAR flow of index $n$ on $\CR$.

We recall the following proposition from \cite{Alev} (see proposition 2.6).

\begin{Proposition}\label{alevCAR}
The CAR flow of rank $n$ on $\CR$ is conjugate to the Clifford flow of rank $2n$.
\end{Proposition}

We point out an error in \cite{ABS} in the following remark.

\begin{Remark}
In section 5, \cite{ABS}, it is claimed that CAR flows of any given
rank are extendable. In fact a `proof' is given, for any $\lambda \in
(0,\frac{1}{2}]$ with $A=\lambda$, that the corresponding
$E_0-$semigroup on $M_A$ is extendable. (When $\lambda\neq
\frac{1}{2}$ they are type III factors.) But we have seen that
Clifford flows are not extendable. This consequently implies, thanks
to proposition \ref{alevCAR}, and the invariance of extendability of
$E_0$-semigroups of $II_1$ factors (where equimodularity - with
respect to the trace - comes for free), that CAR flows on the hyperfinite type II$_1$ factor $\CR$ are not extendable. 

In fact it has been proved in \cite{B} that CAR flows flows on any type $III_\lambda$ factors (considered in \cite{ABS}) are also not extendable.)

\end{Remark}

Let $\Gamma_f(\CH_\C)$ be the full Fock space associated with a Hilbert space $K$. 
For $f \in \CH$, define  $s(f)=\frac{l(f)+l(f)^*}{2}$ 
where $$l(f)\xi=\left\{ \begin{array}{ll} f & \hbox{if}~\xi=\Omega,\\  f\otimes \xi & \hbox{if}~\langle \xi,\Omega\rangle=0. \end{array}\right\}.$$

The von Neumann algebra  $\Phi(\CK)=\{s(f):f\in \CH\}^{\prime\prime} $, 
is isomorphic to the free group factor $L(F_\infty)$ and the vacuum is cyclic and 
separating with $\langle \Omega,x\Omega \rangle=\tau(x)$ (see \cite{voiculescu}) a tracial state on $\Phi(\CK)$.

\begin{Example}
There exists a unique $E_0-$semigroup $\gamma$ on $\Phi(\CK)$ satisfying
 $$\gamma_t(s(f)):=s(S_tf) \qquad (f\in \CH,~t\geq0).$$ This is called the free flow of rank $dim(\CK)$.

Let $\gamma$ be a free flow of any rank. 
It is known - see \cite{Popa} and \cite{Ueda}- that 
$\gamma_t(\Phi(\CK))^\prime \cap \Phi(\CK) = \C 1 $. So it follows from  Theorem 
\ref{rela} that free flow is not extendable.

It is proved in \cite{SM}
that $$ H(t)  = (E^{\gamma_t} \cap E^{\gamma_t^\prime}) = \C S_t.$$
So $H = \{ ( t, \eta ) : \eta \in H(t) \}$ is a
product system. 
This means that free flows provide examples to show that the converse of the Corollary \ref{prodlem} is not true: 
for free flows, the family $\{(t, (E^{\gamma_t}\cap E^{\gamma_t^\prime})~ ) : t\geq 0 \}$ 
forms a product system, but still they are not extendable. 
\end{Example}

\bigskip \noindent
{\bf Acknowledgement:} We thank Kunal Mukherjee for bringing \cite{ABS} to our attention, while 
this work was in progress.
M. I. is supported in part by the Grant-in-Aid for Scientific Research
(B) 22340032, JSPS.  And V.S.S. is happy to thank the J.C. Bose
Fellowship which has helped support this research.

\end{document}